\documentclass[10pt, reqno,twosides]{amsart}

\usepackage{amsmath,amssymb}
\usepackage[english]{babel}
\usepackage[dvips]{graphicx}
\usepackage{indentfirst}
\usepackage{enumerate}
\makeatletter \@addtoreset{equation}{section}

\makeatother \makeatletter

\newtheorem{lemma}{Lemma}[section]

\newtheorem{theorem}[lemma]{Theorem}
\newtheorem{remark}[lemma]{Remark}
\newtheorem{definition}[lemma]{Definition}
\newtheorem{proposition}[lemma]{Proposition}

\newtheorem{assumption}[lemma]{Assumptions}
\newtheorem{properties}[lemma]{Properties}
\newcommand{\R}{\mathbb{R}}
\newcommand{\cc}{\mathbb{C}}
\newcommand{\Tt}{\big( T(t) \big)_{t\geq 0}}
\newcommand{\N}{\mathbb{N}}

\newcommand{\rg}{\hbox{\rm{rg}}}
\newcommand{\Rp}{\rm{Re\,}}

\newcommand{\eps}{\varepsilon}

\newcommand{\spa}{\hbox{\rm{span}}}
\newcommand{\fix}{\hbox{\rm{fix}}}
\newcommand{\co}{\hbox{\rm{co}}}

\newcommand{\ov}{\overline}

\newcommand{\tlim}{\tau \hbox{-} \lim}

\newcommand{\one}{1\!\!\!\;\mathrm{l}}

\renewcommand{\phi}{\varphi}

\newcommand{\sumij}{\sum_{i,j=1}^N}

\newcommand{\Dmax}{D_{\hbox{\small{max}}}}

\newcommand{\loc}{\hbox{\rm{\tiny{loc}}}}

\allowdisplaybreaks

\newcommand{\refer}[5]{{\sc #1}{\ #2}{\em\ #3}{\bf\ #4}{\ #5}}

\newcommand{\refbook}[3]{{\sc #1}{\em\ #2}{\ #3}}

\title{Mean Ergodic Theorems for Bi--continuous Semigroups}
\author{A. A. Albanese}
\address{Dipartimento di Matematica ``Ennio De Giorgi'', Universit\`a del Salento, Via Per Arnesano, C.P. 193, I-73100, Lecce, Italy.}
\email{angela.albanese@unisalento.it}
\author{L. Lorenzi}
\address{Dipartimento di Matematica, Universit\`a degli Studi di Parma, Viale G.P. Usberti 53/A, I-43124 Parma, Italy.}
\email{luca.lorenzi@unipr.it}
\author{V. Manco}
\address{Dipartimento di Matematica ``Ennio De Giorgi'', Universit\`a del Salento, Via Per Arnesano, C.P. 193, I-73100, Lecce, Italy.}
\email{vincenzo.manco@unisalento.it} 

\date{\today}

\begin{document}
\maketitle
\begin{abstract}
\noindent In this paper we study the main properties of the Ces\`aro
means of bi-continuous semigroups, introduced and studied by
K\"{u}hnemund in \cite{Kuhnemund03}. We also give some applications
to Feller semigroups generated by second-order elliptic differential
operators with unbounded coefficients in $C_b(\R^N)$ and to
evolution operators associated with nonautonomous second-order
differential operators in $C_b(\R^N)$ with time-periodic
coefficients.
\end{abstract}

\pagenumbering{arabic}

\section{Introduction}

In the last years, the study of transition Markov semigroups on spaces of bounded continuous  functions or of uniformly continuous bounded functions led to consider a class of semigroups of operators for which the usual strong continuity fails to hold. For instance the Ornstein-Uhlenbeck semigroup on the space of uniformly continuous bounded functions on $\R^N$  or even the  heat semigroup on the space of bounded continuous functions on $\R^N$ are not $C_0$-semigroups with respect to the sup-norm.

It was then natural to look for suitable locally convex topologies weaker than the norm topology to treat the lack of strong continuity.
The results of this paper are given in the general framework introduced by K\"uhnemund  \cite{Kuhnemund03} in this direction. Her approach goes back to
 Cerrai and Priola \cite{Ce1994,Prio99}, who considered, on the space $UC_b(H)$ of bounded uniformly continuous functions on a Hilbert space $H$, the uniform convergence on  compact sets (or the pointwise convergence) of  equibounded sequences (the so-called $\mathcal K$- and $\pi$-convergence).
This type of convergence turned out to be a powerful tool for the study of transition Markov semigroups arising from stochastic differential equations.

In particular, in \cite{Kuhnemund03} K\"uhnemund  gave a unified
approach in an abstract setting, by considering the so-called
bi-continuous semigroups, i.e., semigroups of bounded linear
operators on a Banach space which are locally bi-equicontinuous with
respect to an additional locally convex topology $\tau$ coarser than
the norm topology and such that the orbit maps $t\mapsto T(t)x$ are
continuous with respect to this topology, for any $x\in X$.

In  this paper we introduce and study the main properties and the
convergence of the Ces\`aro means of bi-continuous semigroups. More
precisely, we prove some results on the $\tau$-convergence of the
Ces\`aro means of bi-continuous semigroups with respect to the
topology $\tau$, similar to the classical theorems for
$C_0$-semigroups (see e.g., \cite{Nagel} and \cite{dav}). Moreover,
we apply the results to Feller semigroups generated by autonomous and
non autonomous second-order differential operators with unbounded
coefficients in $C_b(\R^N)$.

\section{Preliminaries}

In this section we introduce some notation and recall definitions
and results from \cite{Kuhnemund03} that we use throughout the
paper.

Let $(X,\|\cdot \|)$ be a Banach space and let $\tau$ be a locally convex topology on $X$. We denote by $X'$ the topological dual of $X$ with respect to the $\|\cdot \|$-topology and by $X'_\tau$ the topological dual of $X$ with respect to $\tau$. We assume also that the following conditions hold.

\begin{assumption} \label{ass:1}
 Let $(X, \|\cdot\|)$ be a Banach space. Suppose that there exists
 a locally convex topology $\tau$ on $X$ such that the following properties are
 satisfied:
\begin{itemize}
 \item every $\|\cdot \|$-bounded $\tau$-Cauchy sequence converges in
 $(X,\tau)$;
\item the topology $\tau$ is coarser than the $\|\cdot \|$-topology;
\item $X'_\tau$ is norming for $X$, i.e., for all $x \in X$
\[
  \|x\| = \sup_{\phi \in X'_\tau, \|\phi\|_{X'}\le 1} |\langle x,\phi\rangle |\,.
\]
\end{itemize}
\end{assumption}

In the following we denote by $P_\tau$ a family of seminorms
defining the locally convex topology $\tau$. Since the topology
$\tau$ is coarser than the $\|\cdot \|$-topology, we may assume
without loss of generality that $p(x)\le \|x\|$ for all $x\in X$ and
$p \in P_\tau$. Following  K\"{u}hnemund  \cite{Kuhnemund03}, we
introduce the following definitions.
\begin{definition}\label{defi:3}
A family $\Tt$ of bounded operators on a Banach space
$(X,\|\cdot\|)$ satisfying Assumptions \ref{ass:1}, is called a
bi-continuous semigroup with respect to $\tau$ if the following
conditions hold:
\begin{enumerate}[\rm (i)]
 \item
 $T(0)=I$ and $T(t+s) = T(t)T(s)$ for all $t,s \ge 0$;
\item
there exist $M\ge 1$ and $\omega \in \R$ such that $\|T(t) x \| \le M e^{\omega t} \|x\|$ for all $t\ge 0$  and $x \in X$;
\item
$\Tt$ is strongly $\tau$-continuous, i.e., the map
  \begin{equation*}
    [0,\infty[\, \ni t \mapsto T(t)x \in X
  \end{equation*}
is $\tau$-continuous for each $x \in X$;
\item
$\Tt$ is locally bi-equicontinuous, i.e., for every $\|\cdot \|$-bounded sequence $(x_n)_{n \in \N}$  which is $\tau$-convergent to $x$ we have
\[
\tau-\lim_{n\to \infty} T(t) x_n = T(t)x,
\]
uniformly on bounded intervals of $[0,\infty[\,$.
\end{enumerate}
\end{definition}

For a bi-continuous semigroup $\Tt$ the number
\[ \omega_0 := \inf \{ \omega \in \R : \hbox{ there exists } M \ge 1\hbox{ such that }  \|T(t)\|_{\mathcal L (X)} \le M e^{\omega t} \hbox{ for all } t \ge 0 \}
\]
is called the {\it growth bound} of the semigroup.

Let $\Tt$ be a bi-continuous semigroup on $(X,\|\cdot\|)$ with
respect to  $\tau$, where $X$ satisfies Assumptions \ref{ass:1}.
Since the space $X_\tau$ is norming for $(X, \|\cdot\|)$, the
$\tau$-Laplace transform $R(\lambda)$ of the semigroup $\Tt$ defined
by
\begin{equation}
R(\lambda)x:=\int_0^\infty e^{-\lambda t}T(t)x\,
dt=\|\cdot\|-\lim_{a\to \infty}\int_0^a e^{-\lambda t} T(t)x \, dt,
\label{star}
\end{equation}
for any $x\in X$, is well defined for all $\lambda\in\cc$ with $\Rp
\lambda>\omega_0$ and $x\in X$ (the integrals have to be understood
as $\tau$-Riemann integrals). Moreover, the family of operators
$(R(\lambda))_{\lambda\in \cc, \Rp \lambda>\omega_0}$ satisfies the
resolvent identity as well as the estimate
 \begin{eqnarray*}
  \|R(\lambda)\|_{\mathcal L (X)} \le \frac{M}{\Rp \lambda -\omega},
\end{eqnarray*}
for all $\lambda \in \cc$ with $\Rp \lambda > \omega>\omega_0$.
Further, since $\lambda R(\lambda)x$ $\tau$-converges to $x$ as
$\lambda\to \infty$, each operator $R(\lambda)$ is injective. Hence
$(R(\lambda))_{\lambda\in\mathbb C,~{\rm Re}\lambda>\omega_0}$ is a
resolvent family. Following the approach in \cite{Kuhnemund03},
  the generator of the bi-continuous semigroup $\Tt$ can thus be defined as the unique closed operator $A:D(A) \subset X \to X$ such that $R(\lambda, A) = R(\lambda)$ for all $\lambda \in \cc$ with $\Rp \lambda > \omega_0$. The operator $(A, D(A))$ is called the $\tau$-generator of $\Tt$.

Another and equivalent way to introduce the $\tau$-generator $(A,
D(A))$ of a bi-continuous semigroup $\Tt$ is to define
\begin{align*}
  D(A) &= \left \{x \in X :\, \sup_{0<t<1} \left \|\frac{T(t)x-x}{t} \right \| <\infty \hbox{ and } \tlim_{t \to 0}\frac{T(t)x-x}{t} \hbox{ exists in } X \right \} \\
Ax &= \tlim_{t \to 0^+}\frac{T(t)x-x}{t}.
\end{align*}
Finally, we recall  the following main properties of $(A,D(A))$ (see \cite[Section 3]{Kuhnemund03}).
\begin{properties} \label{proper_01} Let $(X,\|\cdot\|)$ satisfy Assumptions \ref{ass:1}. Let $\Tt$ be a bi-continuous semigroup on $X$ with respect to $\tau$ and let  $(A,D(A))$ be its $\tau$-generator.
Then, the following properties hold:
\begin{enumerate}[\rm (1)]
 \item
if $ x \in D(A)$, then $T(t)x \in D(A)$ for all $t\ge 0$ and the map
$t\mapsto T(t)x$ is continuously differentiable in $t$ with respect
to  $\tau$ with
\[
\frac{d T (t)}{dt}x = AT(t)x = T(t)Ax,
\]
for all $t\ge 0$;
\item
$x \in D(A)$ and $A x = y$ if and only if $T(t)x-x= \int_0^t T(s)y
\,ds$ for all $t\ge 0$;
\item
$(A,D(A))$ is bi-closed, i.e., for all sequences
$(x_n)_{n\in\N}\subset D(A)$ such that $x_n
\stackrel{\tau}{\rightarrow} x$ and $Ax_n
\stackrel{\tau}{\rightarrow} y$, with $(x_n)_{n\in\N}$ and
$(Ax_n)_{n\in\N}$ $\|\cdot \|$-bounded, we have $x \in D(A)$ and
$Ax=y$;
\item
$D(A)$ is bi-dense in $X$, i.e., for every $x \in X$ there exists a
$\|\cdot \|$-bounded sequence $(x_n)_{n\in\N} \subset D(A)$ which is
$\tau$-convergent to $x$;
\item
the subspace $X_0 = \overline{ D(A)}^{\|\cdot \|}$ is invariant for
the semigroup $\Tt$ and $(T(t)_{|X_0})_{t \ge 0}$ is the strongly
continuous semigroup generated by the part of $A$ in $X_0$.
\end{enumerate}
\end{properties}

For further results on bi-continuous semigroups we refer the reader
to \cite{AlMa-1,AlMa-2,Fa0,Fa,K-2,Kuhnemund03}.

\section{Mean Ergodic Theorems}

Throughout this section, we will always assume that  $(X,\|\cdot\|)$
is a Banach space satisfying Assumptions \ref{ass:1}.

Let $\Tt$ be a bi-continuous contraction semigroup on $(X,\|\cdot\|)$ with respect to $\tau$.
For each $r >0$ we define the operator $C(r):X \to X$   by setting
\begin{eqnarray*}
  C(r) x = \frac{1}{r} \int_0^r T(s) x \,ds,
\end{eqnarray*}
for all $x \in X$, where the integral has to be understood as a $\tau$-Riemann integral. The family of operators $(C(r))_{r>0}$ is called the \emph{Ces\`aro means} of the semigroup $\Tt$.

Here, we study the main properties and the $\tau$-convergence of the
Ces\`aro means $(C(r))_{r>0}$ of the semigroup $\Tt$, starting with
the following results.

\begin{proposition} \label{prop:16} Let $\Tt$ be a bi-continuous contraction semigroup on $(X,\|\cdot\|)$ with respect to $\tau$.
 For every $r>0$ the ergodic mean $C(r)$ is a bounded linear operator on $(X,\|\cdot \|)$. Moreover, for every $x\in X$ and every $\|\cdot\|$-bounded sequence $(x_n)_{n \in \N}$ which is $\tau$-convergent to $x$ we have
\begin{eqnarray*}
 \tlim_{n\to \infty} C(r)x_n = C(r)x,
\end{eqnarray*}
uniformly on bounded intervals of $[0,\infty[\,$.
\end{proposition}
\begin{proof}
The linearity of the operator $C(r)$ for any $r>0$ follows easily
from the linearity of the semigroup. Similarly, the boundedness of
$C(r)$ for any $r>0$ is an easy consequence of the boundedness of
the semigroup and the hypothesis  that $X'_\tau$ is norming for $X$.
Indeed, we have
\begin{equation} \label{eq:61}
 \|C(r)x\|  = \sup_{\phi \in X'_{\tau}, \|\phi\|_{X'}\le 1} |\langle C(r) x, \phi \rangle|
     = \sup_{\phi \in X'_{\tau}, \|\phi\|_{X'}\le 1} \left | \frac{1}{r} \int_0^r \langle T(s) x, \phi \rangle \,ds \right |
 \le \|x\|
\end{equation}
for all $ x \in X$.

Now, fix $x\in X$ and let $(x_n)_{n\in \N}$ be a $\|\cdot
\|$-bounded sequence such that $\tau$-converging to $x$ as
$n\to\infty$. The local bi-equicontinuity of the semigroup implies
that $T(t)x_n$ $\tau$-converges to $T(t)x$ as $n\to\infty$,
uniformly on bounded intervals of $[0,\infty[\,$. Hence, for every
$t_0>0$, $p \in P_\tau$ and $\eps >0$, there exists $\nu \in \N$
such that
\[
 \sup_{t\in [0,t_0]} p(T(t)(x_n-x)) < \eps,
\]
for all $n>\nu$. As a byproduct, we have
\[
 p(C(r)(x_n-x)) = p\left ( \frac{1}{r} \int_0^r T(t)(x_n-x) \,dt \right ) \le \frac{1}{r} \int_0^r p(T(t)(x_n-x)) \, dt
 <\eps,
\]
for every $r \in [0,t_0]$ and $n>\nu$. Thus, the claim is proved.
\end{proof}

\begin{lemma}
 \label{lemma:1}
 Let $\Tt$ be a bi-continuous contraction semigroup on $(X,\|\cdot\|)$ with respect to $\tau$ and let  $(A, D(A))$ be its  $\tau$-generator. Then, we have
\[
 \hbox{\rm ker } A = \hbox{\rm fix } \Tt = \{ x \in X : T(t)x=x \hbox{ for all } t \ge 0\}.
\]
\end{lemma}
\begin{proof}
If $x \in \hbox{ker} A = \{ x \in D(A) : Ax = 0\}$, we have, by Property \ref{proper_01}(1),
\[
 T(t)x - x = \int_0^t T(s)Ax \,ds = 0,
\]
for every $t\ge 0$ and, hence, $x \in \hbox{fix }\Tt$. Conversely,
if $T(t)x=x$ for every $t\ge 0$, then $(T(t)x - x)/t$ clearly
converges to $0$ with respect to $\tau$ and, hence, $ x\in D(A)$
with $Ax=0$.
\end{proof}

\begin{proposition} \label{prop:14} Let $\Tt$ be a bi-continuous contraction semigroup on $(X,\|\cdot\|)$
with respect to $\tau$ and let  $(A, D(A))$ be its $\tau$-generator.
Then,
 the Ces\`aro means $(C(r))_{r> 0}$ of the semigroup $\Tt$  satisfy the following properties:
\begin{enumerate}[\rm (i)]
\item
$C(r)x \in \overline{ \hbox{co} \{ T(t)x : t\ge 0\}}^{\tau}$ for every $x \in X$;
\item
for every $t,r>0$
\begin{equation} \label{eq:62}
(I-T(t)) C(r) = C(r) (I-T(t)) = \frac{1}{r} (I-T(r)) \int_0^t T(s) \, ds;
\end{equation}
\item
if $C(r)x$ $\tau$-converges to some $y\in X$, then $y \in \ker A$.
\end{enumerate}
\end{proposition}
 \begin{proof}
 Properties (i) and (ii) follow directly from the definition of the $\tau$-Riemann integral, the strongly $\tau$-continuity of the semigroup and the semigroup law.

In order to prove (iii), assume that for some $x \in X$ $C(r)x$
$\tau$-converges to $y\in X$ as $r\to\infty$. Let
$(r_n)_{n\in\N}\subset [0,\infty[\,$ be an increasing sequence going
to $\infty$ as $n\to \infty$. Since $\|C(r_n)x\| \le \|x\|$ for each
$n\in \N$ by (\ref{eq:61}), the sequence $(C(r_n)x)_{n\in\N}$ is
$\|\cdot\|$-bounded and $\tau$-converges to $y$. So, from
(\ref{eq:62}) and the local bi-continuity of the semigroup $\Tt$, it
follows that
\begin{eqnarray*}
 (I-T(t)) y = \tlim_{n\to \infty} (I-T(t))C(r_n)x = \tlim_{n \to \infty} \frac{1}{r_n} (I -T(r_n)) \int_0^t T(s)x \, ds =
 0,
\end{eqnarray*}
for every $t >0$. This proves that $ y \in \fix \Tt$.
\end{proof}

\begin{definition} \label{defi_mean_ergodic} A contraction bi-continuous semigroup $\Tt$ on $(X,\|\cdot\|)$ with respect to $\tau$ is said to be $\tau$-mean ergodic if $C(r)x$ $\tau$-converges in $X$ for any $x\in X$.
\end{definition}

For $\tau$-mean ergodic semigroups we can define the operator
\[
 P:X \to X, \quad x \mapsto \tlim_{r \to \infty} C(r)x .
\]
We observe that $P$ is a bounded operator on $X$. Indeed, by (\ref{eq:61}) and by the fact that $X'_\tau$ is norming for $(X,\|\cdot\|)$, we have
\[
 \|Px\|  = \sup_{\phi \in X'_{\tau}, \|\phi\|_{X'}\le 1} |\langle Px, \phi \rangle | = \sup_{\phi \in X'_{\tau}, \|\phi\|_{X'}\le 1} \lim_{r \to \infty} | \langle C(r) x, \phi \rangle|   \le \|x\|,
\]
for all $x \in X$. In the next lemma, we show that $P$ is in fact a
projection. For this reason, $P$ is called the {\it $\tau$-mean
ergodic projection} associated with $\Tt$.

\begin{lemma}
\label{lemma-3.5}
Let $\Tt$ be a bi-continuous contraction semigroup on
$(X,\|\cdot\|)$ with respect to $\tau$ and let $(A,D(A))$ be its
$\tau$-generator. If $\Tt$ is $\tau$-mean ergodic, then the operator
$P$ is a projection on $X$ such that
\begin{equation} \label{eq:65}
 P =T(t)P = PT(t) = P^2,
\end{equation}
for all $ t\ge 0$.
Therefore, $X$ can be decomposed into the direct sum
\[
  X = \rg P \oplus \ker P,
\]
where
\begin{enumerate}[\rm (i)]
 \item
  $\rg P = \fix \Tt = \ker A$;\\[.05mm]
 \item
 $\ov{\spa \{ x-T(t)x : x \in X, t \ge 0\}}^{\|\cdot \|} \subset \ker P\subset  \ov{\spa \{ x-T(t)x : x \in X, t \ge
 0\}}^\tau$;\\[.05mm]
 \item
 $\ov{\spa \{ x-T(t)x : x \in X, t \ge 0\}}^\tau=\ov{\rg A}^\tau$.
\end{enumerate}
\end{lemma}
 \begin{proof}
By Proposition \ref{prop:14}(iii) and Lemma \ref{lemma:1}, it follows that
$Px\in {\rm fix}\Tt$. This clearly shows that $Px=T(t)Px$ for any $t>0$ and any $x\in X$.

To prove that $P$ commutes with $T(t)$ for any $t>0$, it suffices to
observe that $C(n)$ commutes with $T(t)$ for any $n\in\N$ and any
$t>0$ (see \eqref{eq:62}). Since the sequence $(C(n)x)_{n\in\N}$ is
$\|\cdot\|$-bounded and $\tau$-converges to $Px$ as $n\to\infty$,
for any $x\in X$ and the semigroup $\Tt$ is bi-continuous,
$T(t)C(n)x$ $\tau$-converges to $T(t)Px$ as $n\to\infty$ for any
$t>0$. Hence, letting $n\to\infty$ in the relation
$T(t)C(n)x=C(n)T(t)x$, we obtain $T(t)Px=PT(t)x$. Finally, we observe
that
\[
 Px = \frac{1}{r} \int_0^{r} Px\, ds = \frac{1}{r} \int_0^{r} T(s)Px\, ds = C(r)Px
\]
for all $r>0$ and $x\in X$,  thereby implying that $Px=P^2x$ for all $x\in X$ as $r\to \infty$. Thus,  (\ref{eq:65}) is proved.

Let us now prove properties (i)--(iii).

(i). Since $P$ is a bounded projection on $(X, \|\cdot \|)$ we can decompose $X$ into the direct sum
\[
 X = \rg P \oplus \ker P.
\]
As we have already remarked, $Px \in \fix \Tt$ for all $x \in X$. Conversely, let $x \in \fix \Tt$. Then
\[
 C(r)x = \frac{1}{r} \int_0^r T(s)x \, ds = \frac{1}{r} \int_0^r x\, ds = x,
\]
for all $r>0$; letting $r\to \infty$ we obtain $Px=x$ and  hence $x \in \rg P$ so that  (i) is proved.

(ii). Let $y = x - T(t)x$ for some $t>0$ and $x\in X$. By Proposition \ref{prop:14}(ii), we have
\[
 C(r) y = \frac{1}{r} (I-T(r)) \int_0^t T(s)x \, ds;
\]
letting $r \to \infty$, we  deduce  $Py=  0$, i.e.,  $y \in \ker P$,
since $\Tt$ is a contraction semigroup. So, we can conclude that
$\ker P \supset \ov{Z}^{\|\cdot \|}$, where $Z:={\spa}\{ x-T(t)x : x
\in X,\, t \ge 0 \}$.

In order to show the other inclusion, let us fix $\phi \in X'_\tau$
vanishing on $Z$. Then, for all $t\ge 0$ and  $x \in X$ we have
$0=\langle x-T(t)x, \phi\rangle = \langle x, \phi \rangle - \langle
T(t) x, \phi \rangle$  and, therefore,
\[
 \langle C(r)x, \phi \rangle = \frac{1}{r} \int_0^r \langle T(s)x,\phi \rangle \, ds = \langle x,\phi \rangle .
\]
If $ x \in \ker P$, it follows that
\[
 \langle x, \phi \rangle  = \tlim_{r\to \infty }\langle C(r) x, \phi \rangle = \langle Px, \phi \rangle = 0.
\]
We have  so shown that, if any $\phi \in X'_\tau$ vanishes on $Z$,
then it vanishes on $\ker P$, too. Now, suppose that $\ker P$ is not
contained in $\ov Z^\tau$. Then, there exists $x_0 \in \ker P$ such
that $x_0 \notin \ov Z^\tau$. By the Hahn-Banach theorem, there
exists $\phi \in X'_\tau$ such that
\begin{equation} \label{eq:69}
 \langle x_0, \phi \rangle = 1 \quad \hbox{and} \quad \langle x-T(t)x, \phi \rangle = 0,
\end{equation}
for all $x \in X$ and $t\ge 0$; but (\ref{eq:69}) is a contradiction, because, by the previous argument, $\phi(x_0)$ should be $0$.

(iii). Let us fix $y\in \rg A$ and let $x\in D(A)$ be such that
$y=Ax$. As $(A,D(A))$ is the $\tau$-generator of $\Tt$, we have
\[
y=\tau-\lim_{t\to 0^+}\frac{T(t)x-x}{t}
\]
and so $y\in \ov{Z}^\tau$. Thus, we can conclude that $\ov{\rg
A}^\tau\subset \ov{Z}^\tau$. To show the converse inclusion, take
any $y=x-T(t)x\in Z$. As $D(A)$ is bi-dense in $X$ by Property
\ref{proper_01}(4), there exists a $\|\cdot\|$-bounded sequence
$(x_n)_{n\in\N}\subset D(A)$ which is $\tau$-convergent to $x$,
thereby implying that
\begin{equation}\label{eq:nuova}
y=\tau-\lim_{n\to \infty}(x_n-T(t)x_n),
\end{equation}
as $\Tt$ is a bi-continuous semigroup with respect to $\tau$. On the
other hand, by Property \ref{proper_01}(2)
\begin{equation}\label{eq:nnuova}
x_n-T(t)x_n=\int_0^tAT(s)x_n\, ds\in \ov{\rg A}^\tau,
\end{equation}
for all $n\in\N$. Combining (\ref{eq:nuova}) and (\ref{eq:nnuova}), we obtain that $y\in \ov{\rg A}^\tau$. So, $\ov{Z}^\tau\subset
\ov{\rg A}^\tau$.
\end{proof}

For strongly continuous semigroups it is well known that the
ergodicity of the semigroup can be stated equivalently in terms of
the strong convergence of the operator $\lambda R(\lambda,A)$ as
$\lambda\to 0^+$ (see e.g., \cite[Theorem 5.1]{dav85}). We are going
to show that an analogous property is satisfied by bi-continuous
semigroups.
To achieve this goal, we begin by proving the following theorem,
which is the proper version in $]0,\infty[\,$ of Wiener's theorem
(see \cite[Theorem 3 p. 357]{yos} or  \cite[Theorem 9.7]{rud87}).

\begin{theorem}[Wiener] \label{theo:11}
Let $f \in L^1(]0,\infty[)$ be such that
  \begin{equation}\label{eq:72}
    \int_0^\infty f(x) x^{-i\xi}\, dx \ne 0,
  \end{equation}
for every $\xi \in \R$. Then, the linear span of the set $\{f_\alpha\}_{\alpha >0}$ of functions defined by $f_\alpha (x) = f(\alpha x)$ is dense in $L^1(]0,\infty[)$.
\end{theorem}
\begin{proof}
 Consider the map $F: L^1(\R) \to L^1(]0,\infty[)$ defined by
\begin{equation*}
  Fg(x)= \frac{1}{x}g(\log x),\qquad  x >0,
\end{equation*}
for all $g \in L^1(\R)$. Clearly, $F$ is an isometry with inverse
$F^{-1}h(y)= e^y h(e^y)$ for every $y \in \R$ and any $h \in
L^1(]0,\infty[)$. Since $f\in L^1(]0,\infty[)$ satisfies
(\ref{eq:72})
\begin{align}\label{eq:tr}
   \int_{\R} \mathcal (F^{-1}f) (y) e^{-i\xi y}\, dy = \int_{\R} e^y f(e^y) e^{-i\xi y} \, dy = \int_0^\infty f(x) x^{-i\xi} \, dx \ne 0,
\end{align}
for every $\xi \in \R$, i.e., the Fourier transform of $F^{-1}f$ nowhere vanishes in $\R$.

Now, let us fix $h \in L^1(]0,\infty[)$ and set $g=  F^{-1}h\in L^1(\R)$. Since $F^{-1}f$ satisfies (\ref{eq:tr}), we can apply  Wiener's theorem, obtaining that for every $\eps >0$ there exist $n\in\N$, $\beta_1, \dots, \beta_n \in \cc$ and $\sigma_1, \dots, \sigma_n \in \R$ such that
\begin{equation} \label{eq:73}
   \int_\R \left | g (y) - \sum_{i=1}^n \beta_i ( F^{-1} f)(y+\sigma_i)\right | \, dy < \eps .
\end{equation}
Setting $\alpha_i=e^{\sigma_i}$ and $\gamma_i=\beta_i \alpha_i$ for each $i=1,\ldots, n$,  (\ref{eq:73}) gives
\begin{eqnarray*}
 \int_0^\infty \left | h(x) - \sum_{i=1}^n \gamma_i f(\alpha_i x) \right | \, dx < \eps .
\end{eqnarray*}
This completes the proof.
\end{proof}

 We can now show the following fact.

\begin{lemma}\label{lem:17}
Let $K \in L^1(]0,\infty[)$ be such that
\begin{eqnarray*}
 \int_0^\infty K(y)\,  y^{-i \xi}\, dy \ne 0,
\end{eqnarray*}
for all $\xi \in \R$. If $\Phi\colon [0,\infty[\to X$ is a
$\|\cdot\|$-bounded, $\tau$-continuous function such that
\begin{equation} \label{eq:74}
\tlim_{\lambda \to 0^+} \lambda \int_0^\infty K(\lambda t) \Phi(t)
\, dt = a \int_0^\infty K(t) \, dt,
\end{equation}
for some $a\in X$, then
\begin{eqnarray*}
  \tlim_{\lambda \to 0^+} \lambda \int_0^\infty f(\lambda t) \Phi(t) \, dt = a \int_0^\infty f(t) \, dt
\end{eqnarray*}
for all $f \in L^1(]0,\infty[)$.
\end{lemma}
 \begin{proof}
 We first prove  that property (\ref{eq:74}) holds for every  $f\in \spa \{K_\alpha : \alpha>0\}\subset L^1(]0,\infty[)$, where $K_\alpha(x)=K(\alpha x)$, for all $\alpha, x>0$. For this purpose, we fix $\alpha>0$ and $p \in P_\tau$ and observe that
\begin{equation*}
 p\left ( \lambda \int_0^\infty K_\alpha(\lambda t) \Phi(t)\, dt - a \int_0^\infty K_\alpha(t)\, dt \right ) = \frac{1}{\alpha} p\left (\mu \int_0^\infty K(\mu t) \Phi(t) \, dt - a \int_0^\infty K(s)\, ds\right ),
\end{equation*}
where we have set $\mu = \lambda \alpha$ in the first integral and $s=\alpha t$ in the second one. So, letting $\lambda \to 0^+$ (and hence $\mu \to 0^+$), we obtain
\begin{equation*}
  p\left ( \lambda \int_0^\infty K_\alpha (\lambda t) \Phi(t)\, dt - a \int_0^\infty K_\alpha(t)\, dt \right ) \to 0.
\end{equation*}
By the arbitrariness of $p \in P_\tau$ and by linearity, property (\ref{eq:74}) holds for every $f \in \spa\{K_\alpha : \alpha >0\}$.

Now, let $f \in L^1(]0,\infty[)$. By Lemma \ref{theo:11} there exists a sequence $(\sigma_n)_{n\in\N} \subset \spa \{K_\alpha:\alpha>0\}$ which converges to $f$ in $L^1(]0,\infty[)$. Then, for every $p \in P_\tau$ we have
\begin{align}
  p\left (\lambda \int_0^\infty f(\lambda t) \Phi(t) \, dt\right . & \left . - a \int_0^\infty f(t) \, dt \right )\nonumber \\ & \le  p \left ( \lambda \int_0^\infty[ f(\lambda t)-\sigma_n(\lambda t)] \Phi(t) \, dt \right )\label{eq:75} \\  &\quad  + p \left ( \lambda \int_0^\infty \sigma_n(\lambda t) \Phi(t) \, dt- a \int_0^\infty \sigma_n(t)\, dt \right )\label{eq:76} \\ & \quad + p\left (a\int_0^\infty[ \sigma_n(t) - f(t)] \, dt \right ). \label{eq:77}
\end{align}
Observe that (\ref{eq:76}) tends to zero as $\lambda \to 0^+$ for
all $n\in \N$ as $(\sigma_n)_{n\in \N}\subset \spa\{K_\alpha :
\alpha>0\}$. On the other hand, (\ref{eq:75}) and (\ref{eq:77}) tend
to zero as $n \to \infty$ by the dominated convergence theorem
because $\sigma_n \to f$ in $L^1(]0,\infty[)$ and $\Phi$ is
$\|\cdot \|$-bounded. It is now immediate to conclude that the first
side of the previous chain of inequalities vanishes as $\lambda\to
0^+$, accomplishing the proof.
\end{proof}

We are now able to show that the $\tau$-mean ergodicity of a bi-continuous contraction semigroup $\Tt$ can be also characterized in terms of the $\tau$-convergence of its resolvent operator $R(\lambda,A)$.

\begin{proposition} \label{prop:15}
 Let $\Tt$ be a contraction semigroup on $(X,\|\cdot\|)$, bi-continuous with respect to $\tau$, and let $(A,D(A))$ be its $\tau$-generator. Then, the  following statements are equivalent:
\begin{enumerate}[\rm (i)]
 \item
$\Tt$ is $\tau$-mean ergodic, i.e.,  $C(r)x$ $\tau$-converges in $X$
as $r\to\infty$, for every $x \in X$;
\item
$\lambda R(\lambda,A)x$ $\tau$-converges in $X$ as $\lambda\to 0^+$,
for every $x \in X$.
\end{enumerate}
Moreover, if one of these limits exists, then $\tlim_{r \to \infty} C(r)x = \tlim_{\lambda \to 0^+} \lambda R(\lambda,A)x$; consequently, $Px = \tlim_{\lambda \to 0^+} \lambda R(\lambda,A)x$.
\end{proposition}
 \begin{proof}
 (i)$\Rightarrow$(ii). Assume that $\Tt$ is $\tau$-mean ergodic and fix $x \in X$. By Definition \ref{defi_mean_ergodic}, $\tlim_{r \to \infty} C(r)x = Px \in X$. Since the semigroup is contractive, it has growth bound $\omega_0\le 0$ and, hence, $\lambda \in \rho(A)$ for every $\lambda >0$ with
\begin{equation} \label{eq:63}
 R(\lambda,A)x = \int_0^\infty e^{-\lambda s} T(s) x \, ds,
\end{equation}
for any $x\in X$. Integrating by parts in (\ref{eq:63}) we obtain
\begin{align*}
 \lambda R(\lambda,A) x &=  \lambda^2 \int_0^\infty s e^{-\lambda s} C(s)x \, ds.
\end{align*}
As $\lambda^2 \int_0^\infty s e^{-\lambda s} \: ds = 1$, for any $p \in P_\tau$ we have
\begin{align*}
p( \lambda R(\lambda, A) x -Px) & = p \left (\lambda^2 \int_0^\infty
s e^{-\lambda s} C(s)x \, ds - Px \right )\\ &=
    p \left (\lambda^2 \int_0^\infty s e^{-\lambda s} [ C(s)x -Px] \, ds  \right ) \\ &=
    p \left ( \int_0^\infty t e^{-t} [ C(t/\lambda)x -Px] \, dt  \right ) \\ &\le
    \int_0^\infty t e^{-t} p(C(t/\lambda)x -Px) \, dt ,
\end{align*}
where by (i) $p(C(t/\lambda)x -Px) \to 0$ as $\lambda \to 0^+$ for every $t>0$. By the dominated convergence theorem, it follows that $p(\lambda R(\lambda, A)x -Px)\to 0$ as $\lambda \to 0^+$. This completes the proof.

(ii)$\Rightarrow$(i). Fix $x\in X$. By assumption $\lambda
R(\lambda,A)x$ $\tau$-converges to some $a\in X$ as $\lambda\to
0^+$. Let us observe that, by setting $\lambda = 1/r$, we have
\begin{equation} \label{eq:79}
 C(r)x = \lambda \int_0^{1/\lambda} T(s)x \, ds = \lambda \int_0^\infty \chi_{[0,1]}(\lambda s) T(s) x \, ds.
\end{equation}
So, it suffices to prove that
\begin{equation*}
  \tlim_{\lambda \to 0^+} \lambda \int_0^\infty \chi_{[0,1]}(\lambda s) T(s)x \, ds= a.
\end{equation*}
In order to do this, we observe that for every $\lambda >0$
\begin{eqnarray*}
  \int_0^\infty e^{-\lambda s} s^{-i\xi} \, ds = \lambda^{i\xi -1} \int_0^\infty t^{-i\xi} e^{-t} \, dt = \lambda^{i\xi -1} \Gamma(1-i\xi) \ne 0, \quad \xi \in \R,
\end{eqnarray*}
because the Euler function $\Gamma$ nowhere vanishes  in $\cc \setminus \mathbb Z^-$ (for the integral representation of the Euler function and its main properties we refer to \cite[Chapter 2, Sections 8 and 15]{Rain}). By assumption
\begin{equation*}
  \tlim_{\lambda \to 0^+} \lambda \int_0^\infty e^{-\lambda s} T(s)x \, ds = a = a \int_0^\infty e^{-s} \, ds.
\end{equation*}
Since $T(\cdot)x$ is $\|\cdot\|$-bounded and $\tau$-continuous, and
$f(t)=e^{-t}\in L^1(]0,\infty[)$, we can then apply Lemma
\ref{lem:17} to conclude that there exists
 \begin{equation*}
   \tlim_{\lambda \to 0^+} \lambda \int_0^\infty \chi_{[0,1]}(\lambda s) T(s)x \, ds = a \int_0^\infty \chi_{[0,1]}(s) \, ds = a.
 \end{equation*}
By (\ref{eq:79}), this completes the proof.
\end{proof}

\begin{remark}
{\rm  Note that also the implication (ii)$\Rightarrow$(i) can be proved applying Wiener's theorem. We have preferred to
provide here a direct proof for the sake of simplicity.}
\end{remark}

The $\tau$-mean ergodicity of a bi-continuous contraction semigroup
can be also described by the following series of properties.

\begin{theorem} \label{theo:1}
Let $\Tt$ be a contraction semigroup on  $(X,\|\cdot\|)$,
bi-continuous with respect to $\tau$, and let $(A,D(A))$ be its
$\tau$-generator. Consider the following assertions:
\begin{enumerate}[\rm (a)]
 \item
 $\Tt$ is $\tau$-mean ergodic;
\item
the Ces\`aro means $(C(r))_{r >0}$ converges in the $\tau$-weak
operator topology as $r \to \infty$, i.e., $\langle C(r)x, \phi
\rangle$ converges as $r\to \infty$ for every $\phi \in X'_\tau$ and
$x\in X$;
\item
for every $x \in X$ there exists an unbounded sequence $(r_n)_{n\in
\N} \subset \R^+$ such that $(C(r_n)x)_{n\in \N}$ has a $\tau$-weak
accumulation point in $X$;
\item
for every $x \in X$ one has $\ov{ \co\{T(t)x : t\ge 0\}}^\tau \cap
\fix \Tt \neq \varnothing$;
\item
the fixed space $\fix \Tt = \ker A$ separates the dual fixed space
$\fix (T(t)')_{t \ge 0} = \ker A'$ in $X'_\tau$.
\end{enumerate}
Then, {\rm (a)}$\Rightarrow${\rm (b)}$\Rightarrow${\rm (c)} and {\rm
(a)}$\Rightarrow${\rm (d)}$\Rightarrow${\rm (e)}. Further, if for
every $\|\cdot \|$-bounded sequence $(x_n)_{n\in \N}$, which is
$\tau$-weak convergent to $x$, it holds that
\begin{equation} \label{eq:70}
 \sigma(X,X'_\tau) \hbox{-}\lim_{n\to\infty} T(t)x_n = T(t)x ,
\end{equation}
then {\rm (c)}$\Rightarrow${\rm (d)}. Finally, if the semigroup is equicontinuous, i.e.,
for every $p \in P_\tau$ there exist $c_p\ge 1$ and $ q \in P_\tau$ so that
\begin{equation} \label{eq:81}
  p(T(t)x) \le c_p q(x),
\end{equation}
for all $x \in X$, $t\ge 0$, then {\rm (e)}$\Rightarrow${\rm (a)}.
\end{theorem}
 \begin{proof}
The proof of (a)$\Rightarrow$(b)$\Rightarrow$(c) is immediate and, hence, we skip it.

(c)$\Rightarrow$(d). Let $x \in X$ and $(r_n)_{n \in \N}\subset \R^+$ be an unbounded sequence as in (c).
Then, up to a subsequence we can assume that $C(r_n)x$ $\tau$-weak
converges to some $y\in X$ as $n\to\infty$. By condition
(\ref{eq:70}), for any $t>0$, $(I-T(t))C(r_n)x$ converges to $y-
T(t)y$ as $n\to \infty$ with respect to the $\tau$-weak topology.
Moreover, by Proposition \ref{prop:14}(ii) we have
\[
 \| (I-T(t)) C(r_n)x \| \le \frac{2t}{r_n} \|x\|,
\]
for all $n \in \N$. Hence, $(I-T(t))C(r_n)x \to 0$  in $(X,\|\cdot\|)$ as $n\to \infty$. By the uniqueness of the limit, we obtain that $y-T(t)y = 0$. The arbitrariness of $t\ge 0$ implies that $y\in \fix \Tt$.

On the other hand, by Proposition \ref{prop:14}(i), $C(r_n)x \in \overline{ \hbox{co} \{ T(t)x : t \ge 0 \} }^\tau $ for all $n\in \N$, which is clearly a $\tau$-weak closed set of $X$, therefore implying that $y \in \ov{ \hbox{co} \{T(t)x : t \ge 0 \} }^\tau$. So, (c)$\Rightarrow$(d) is proved.

(a)$\Rightarrow$(d). The proof is along the lines of the proof of
the implication (c)$\Rightarrow$(d). Hence, we skip the details.

(d)$\Rightarrow$(e). Let us fix $x', y' \in \ker A' = \fix
(T(t)')_{t \ge 0}$ such that $x',y' \in X'_\tau$ and $x' \neq y'$.
Then, there exists $x_0 \in X$ such that
\[
 \langle x_0, x'\rangle \neq \langle x_0, y'\rangle .
\]
By assumption there exists $\ov x \in \ov{\hbox{co} \{T(t)x_0 : t \ge 0\}}^\tau \cap \fix \Tt$. Then, for every $y \in \hbox{co} \{T(t)x_0 : t \ge 0\}$, i.e., for $y = \sum_{i=1}^k\lambda_i T(t_i) x_0$ with $\sum_{i=1}^k \lambda_i=1$, we have
\begin{align*}
 \langle y,x' \rangle   = \left\langle \sum_{i=1}^k \lambda_i T(t_i)x_0, x' \right\rangle = \sum_{i=1}^k \lambda_i \langle x_0, T(t_i)'x' \rangle    = \sum_{i=1}^k \lambda_i \langle x_0,x' \rangle = \langle x_0,x' \rangle
\end{align*}
and, analogously, we have
\[
 \langle y,y'\rangle = \langle x_0,y' \rangle .
\]
Since $x',y'$ are $\tau$-continuous in $X$, it follows that
$\langle y,x'\rangle = \langle x_0,x' \rangle$ and $\langle y,y'\rangle = \langle x_0,y' \rangle$
for all $y \in \ov{\hbox{co} \{T(t)x_0 : t \ge 0\}}^\tau$. In particular, since $\ov x \in \ov{\hbox{co} \{T(t)x_0 : t \ge 0\}}^\tau$, we have
\[
 \langle \ov x, x' \rangle = \langle x_0, x' \rangle \neq \langle x_0, y' \rangle = \langle \ov x, y' \rangle ,
\]
i.e., $\fix \Tt$ separates $\fix (T(t)')_{t\ge 0}$ in $X'_\tau$.

(e)$\Rightarrow$(a). We first observe  that condition (\ref{eq:81}) implies that
\begin{equation*}
  p(C(r)x) \le c_p q(x),
\end{equation*}
for all $x\in X$, $r>0$. Let us consider the subspace
\[
 G = \fix \Tt \oplus \spa \{ x-T(t)x : x\in X, t \ge 0\}
\]
of $X$. (To show that the two subspaces constituting $G$ intersect
in $0$ only, it suffices to recall that $\fix\Tt=\rg P$, ${\rm
ker}\,P\supset \spa \{ x-T(t)x : x\in X, t \ge 0\}$ and $\ker P\cap
\rg P=\{0\}$)

Fix $x' \in X'_\tau$ vanishing on $G$. Then, for every $t\ge 0$, $x \in X$, we have
\[
 \langle x-T(t)x,x' \rangle = 0 \Longleftrightarrow \langle x , x' - T'(t)x' \rangle = 0,
\]
i.e., $x' \in X'_\tau \cap \fix (T'(t))_{t\ge 0}$. Since $\langle x, x'\rangle= 0$ for all $x \in \fix \Tt$ too, and by assumption $\fix \Tt$ separates the space $X'_\tau \cap \fix (T'(t))_{t \ge 0}$, we can then conclude that $x'=0$ (otherwise there would exist $x_0 \in \fix \Tt$ such that $\langle x_0, x'\rangle \ne 0$). By the arbitrariness of $x' \in X'_\tau$, we obtain that $\ov{G}^\tau =X$.

To prove that $\Tt$ is $\tau$-mean ergodic, it now suffices to prove that
\begin{equation} \label{eq:82}
C(n)x = \frac{1}{n} \int_0^n T(t)x \,dt\to y,
\end{equation}
in $X_\tau$ as $n\to \infty$ for any $x\in X$. Indeed, if (\ref{eq:82}) holds true, then for every $r>0$
\begin{equation*}
  C(r)x = \frac{1}{r} \int_0^n T(t)x \, dt + \frac{1}{r} \int_n^{n+\alpha}T(t)x \, dt ,
\end{equation*}
where $n=[r]$ and $\alpha = r-n \in [0,1[\,$. Therefore,
\begin{align}\label{eq:83}
  C(r)x &= \frac{n}{r} C(n)x + \frac{1}{r} \int_0^\alpha T(t+n)x \, dt \nonumber \\
  &= \left ( 1-\frac{\alpha}{r} \right ) C(n)x + \frac{\alpha}{r} T(n) C(\alpha)x \nonumber \\
  &= \left ( 1-\frac{\alpha}{r} \right ) \left [C(n)x + \frac{\alpha}{n} T(n) C(\alpha)x \right ].
\end{align}
Since $1-\alpha/r \to 1$ as $r \to \infty$, $C(n)x \to y$ in $X_\tau$ as $n\to \infty$ and $\|\alpha C(\alpha)T(n) x\| \le \|x\|$, letting $r\to \infty$ (hence $n=[r] \to \infty$) we then obtain that $C(r)x \to y$ in $X_\tau$.

To prove (\ref{eq:82}), let us fix $p \in P_\tau$ and $\eps >0$, and
let $\ov x \in G$ satisfy $q(x-\ov x)< \eps c_p /3$ ($c_p$ and $q$
are chosen according to (\ref{eq:81})). Since $\ov x\in G$,  the
$\|\cdot\|$-bounded sequence $(C(n)\ov x)_{n\in\N}$ converges in
$(X, \|\cdot \|)$, say to $ \ov y$, and hence $(C(n)\ov x)_{n\in\N}$
is a Cauchy sequence in $(X, \|\cdot \|)$. So, there exists $n_0 \in
\N$ such that
\begin{equation*}
    \| C(n)\ov x - C(m) \ov x\| < \frac{\eps}{3}
\end{equation*}
for all $n, m \ge n_0$. Combining all these facts, we get that, for all $m, n \ge n_0$,
\begin{align*}
p(C(n)x - C(m)x) &\le p(C(n)x - C(n)\ov x) + p(C(n)\ov x - C(m)\ov x) + p(C(m)\ov x-C(m)x)\\
& \le 2 c_p q(x-\ov x) + \|C(n)\ov x - C(m) \ov x\| \\ &\le 2 c_p \frac{\eps}{3c_p} + \frac{\eps}{3} < \eps.
\end{align*}
The arbitrariness of $\eps$ and $p$ imply that $(C(n))_{n\in \N}$ is a $\|\cdot\|$-bounded $\tau$-Cauchy sequence and, hence, it converges in $X_\tau$ by Assumptions \ref{ass:1}.
\end{proof}

\begin{remark}{\rm
Assume that the following conditions are satisfied:
\begin{enumerate}[\rm (1)]
\item
there exists another norm $|\cdot |$ on $X$ such that the inclusion $(X,\tau)\hookrightarrow (X,|\cdot |)$ is continuous;
\item
the topology $\tau$ and the  $|\cdot |$-topology coincide on the bounded sequences of $(X, \|\cdot\|)$;
\item
there exists $M\ge 1$ such that $|T(t)x|\le M |x|$ for all $t \ge 0$ and $x \in X$.
\end{enumerate}
Then the implication (e)$\Rightarrow$(a) holds even without the assumption on the boundedness of $\Tt$ with respect to $\tau$.
Indeed, as in the proof of Theorem \ref{theo:1}, it suffices to show only that $C(n)x \to y$ in $X_\tau$ (for some $y \in X$).
Fix $\eps >0$ and let $ \ov x \in G$ satisfy $|x-\ov x| < \eps / 3M$ (note that $\ov G^\tau = X$ implies $\ov G^{|\cdot |} = X$ by (1)). Since $(C(n)\ov x)_{n\in \N}$ is a Cauchy sequence in $(X, \|\cdot \|)$ and hence in $(X, |\cdot |)$, there exists $n_0 \in \N$ such that
\begin{equation*}
 |C(n)\ov x - C(m) \ov x| \le c\|C(n)\ov x  - C(m)\ov x\| < \frac{\eps}{3}
\end{equation*}
for all $n, m \ge n_0$ (note that by Assumptions \ref{ass:1} and
condition (1) above  there exists a constant $c>0$ such that $|x|\le
c\|x\|$ for every $x\in X$). Further, notice that $|C(r)x|\le M|x|$
for any $x\in X$ and any $r>0$. Indeed, the function $t\mapsto
T(t)x$ is continuous with respect to $|\cdot |$ in $[0,\infty[$\,.
Hence, it is $|\cdot|$-Riemann integrable in any bounded interval of
$[0,\infty[\,$. Moreover, the Riemann sums whose limit defines
$C(r)$ constitute a $\|\cdot\|$-bounded set in $X$. Since $\tau$ and
$|\cdot|$ define the same topology on the $\|\cdot\|$-bounded sets
of $X$ by (2),
\begin{eqnarray*}
\tau-\int_0^tT(s)xds=|\cdot|-\int_0^tT(s)xds,
\end{eqnarray*}
for any $t>0$ and any $x\in X$. Hence,
\begin{eqnarray*}
|C(r)x|=\frac{1}{r}\left |\int_0^rT(s)xds\right |\le M|x|,
\end{eqnarray*}
for any $r>0$ and any $x\in X$, where $M$ is the constant in (3).

Taking all the previous remarks into account, we can infer that
\begin{align*}
 |C(n)x- C(m)x| &\le |C(n)\ov x - C(n)x|+|C(n) \ov x - C(m) \ov x|+|C(m) \ov x - C(m)x| \\
        &\le 2M |x- \ov x| + \frac{\eps}{3} < \eps.
\end{align*}
This means that $(C(n)x)_{n \in \N}$ is a $|\cdot |$-Cauchy sequence
and, hence, a $\tau$-Cauchy sequence by (2), as the sequence
$(C(n)x)_{n \in \N}$ is $\|\cdot\|$-bounded. Assumptions \ref{ass:1}
imply that $C(n)x$ $\tau$-converges in $X$. }
\end{remark}

Denote by $\tau_m$ the finest locally convex topology on $X$
agreeing with $\tau$ on the bounded sets of $(X, \|\cdot\|)$.
Actually, $\tau_m$ is the so-called mixed topology and we refer to
\cite{wiw} for more details. Then, Theorem \ref{theo:1} can be
reformulated as follows.

\begin{theorem} Let $\Tt$ be a contraction semigroup on $(X,\|\cdot\|)$ which is bi-continuous with respect to $\tau$ and let $(A, D(A))$ be its $\tau$-generator. Consider the following properties:
\begin{enumerate}[\rm (a)]
 \item
 $\Tt$ is $\tau$-mean ergodic;
\item
the Ces\'aro means $(C(r))_{r>0}$ converge in the $\tau_m$-weak
operator topology as $r \to \infty$, i.e., $\langle C(r)x, \phi
\rangle$ converges as $r\to \infty$ for every $\phi \in X'_{\tau_m}$
and $x \in X$;
\item
for every $x \in X$ there exists an unbounded sequence $(r_n)_{n\in \N} \subset \R^+$ such that $(C(r_n)x)_{n\in \N}$ has a $\tau_m$-weak accumulation point in $X$;
\item
for every $x \in X$ one has $\ov{ \co\{T(t)x : t\ge 0\}}^{\tau_m} \cap \fix \Tt \neq \varnothing$;
\item
the fixed space $\fix \Tt = \ker A$ separates the dual fixed space $\fix (T(t)')_{t \ge 0} = \ker A'$ in $X'_{\tau_m}$.
\end{enumerate}
Then, {\rm (a)}$\Rightarrow${\rm (b)}$\Rightarrow${\rm (c)}, {\rm
(d)}$\Rightarrow${\rm (e)}. Moreover, if the topology $\tau$ is
metrizable, then {\rm (c)}$\Rightarrow${\rm (d)}
\end{theorem}

\begin{proof}
It suffices to recall that the topologies $\tau$ and $\tau_m$
coincide on the $\|\cdot\|$-bounded sets. This assures that
\begin{center}
 $C(r)x$ $\tau$-converges $\Longleftrightarrow$ $C(r)x$ $\tau_m$-converges
\end{center}
as $r\to \infty$ via identity (\ref{eq:83}). Therefore, (a) implies that $\Tt$ is $\tau_m$-mean ergodic. Then the $\tau_m$-mean ergodicity clearly implies (b).

Next, to prove (b)$\Rightarrow$(c) and (d)$\Rightarrow$(e) it suffices to proceed as in Theorem \ref{theo:1}.\\
(c)$\Rightarrow$(d). Since $\Tt$ is locally $\tau$-bi-equicontinuous
and $\tau$ is metrizable, we can apply  \cite[Corollary 2.2.4]{wiw}
(see also \cite[2.2.2]{wiw}) to conclude  that $T(t): X_{\tau_m} \to
X_{\tau_m}$ is continuous and hence,
$\sigma(X,X'_{\tau_m})-\sigma(X,X'_{\tau_m})$
continuous. Therefore, it suffices again to proceed as in the proof
of the corresponding implication in Theorem \ref{theo:1}.
\end{proof}

\section{Examples}

In this section, we show some application of the results in the
previous sections to semigroups associated with second-order
elliptic operators $A$ with possibly unbounded coefficients. We
consider both the case of autonomous and nonautonomous operators.

\subsection{The autonomous case}

In this subsection, we will consider semigroups of linear continuous
operators on the Banach space $(C_b(\R^N),\|\cdot \|_\infty)$ that
are bi-continuous with respect to the  topology $\tau_c$ of the
uniform convergence on compact sets of $\R^N$.

Let $A$ be the second-order elliptic partial differential operator (with possibly unbounded coefficients) defined on smooth functions $u$ by
\begin{equation}\label{eq:1}
    Au(x) = \sumij q_{ij}(x)D_{ij} u(x) + \sum_{i=1}^N b_i(x) D_i u(x), \qquad x \in \R^N.
\end{equation}
Autonomous differential operators with unbounded coefficients on
$\R^N$ have been investigated intensively in recent years after the
seminal paper \cite{DaLu} (see, e.g.,
\cite{AlLoMa,AlMa-3,AlMa-4,libro,Ce1994,lun98,MPW2002,MPW2002-1,Prio99} and the
references therein) both in spaces of continuous functions and in
suitable $L^p$-spaces. We now recall some results on the
differential operator $A$ defined  in (\ref{eq:1}), that we need in
what follows, from the survey paper \cite{MPW2002}, under the
following (minimal) assumptions on the coefficients:
\begin{enumerate}[1)]
\item
the coefficients $q_{ij}, b_i,c \in C^{0,\alpha}_{\loc} (\R^N)$ for
some $0 <\alpha <1$;
\item
$q_{ij}=q_{ji}$ for every $i,j =1, \dots , N$ and the ellipticity
condition
\begin{eqnarray*}
\sumij q_{ij}(x) \xi_i\xi_j \ge \nu (x) |\xi|^2, \qquad x,\ \xi \in
\R^N,
\end{eqnarray*}
is satisfied, where $\inf_{K}\nu >0$ for every compact set $K\subset \R^N$.
\end{enumerate}

In \cite{MPW2002} the authors presented a very elegant construction of a semigroup $(T(t))_{t\geq 0}$ of positive contractions on $C_b(\R^N)$ associated with the operator $A$. In fact, for any $f\in C_b(\R^N)$, $T(t)f$ gives, for positive $f\in C_b(\R^N)$, the value at $t>0$ of the minimal positive bounded classical solution of the parabolic problem
\[
\left\{ \begin{array}{ll}
D_tu(t,x)=Au(t,x), & \mbox{ $t>0$, $x\in\R^N$,}\\[1mm]
u(0,x)=f(x), &\mbox{ $x\in\R^N$.}
\end{array}\right.
\]
Here, by classical solution, we mean a function $u\in
C^{1,2}(]0,\infty[\,\times\R^N)\cap C([0,\infty[\times\R^N)$ which
solves the Cauchy problem. Moreover,
\begin{equation}
T(t)f(x)=\int_{\R^N}p(t,x,y)f(y)\,dy,\qquad t>0,\,\,x\in\R^N,
\label{repres}
\end{equation}
for any $f\in C_b(\R^N)$, where $p$ is a positive and smooth
function (for further properties of $p$ we refer to \cite[Theorems
4.4, 4.5]{MPW2002}). Using this representation formula, it is
immediate to check that the semigroup $(T(t))_{t\geq 0}$ is
irreducible, that is $T(t)f(x)>0$ for every $t>0$, $x\in\R^N$
whenever $f\geq 0$, $f\not =0$. Moreover, \eqref{repres} combined
with Schauder interior estimates show that $(T(t))_{t\geq 0}$ has
the strong Feller property, i.e., $T(t)f\in C_b(\R^N)$ for every
bounded Borel function $f$. In particular, $(T(t))_{t\geq 0}$ is not
strongly continuous in $C_b(\R^N)$, but $T(t)f$ tends to $f$ as
$t\to 0$, uniformly on compact subsets of $\R^N$ (this is a typical
behaviour for semigroups associated with elliptic operators with
unbounded coefficients). This property assures that the semigroup
$\Tt$ is bi-continuous with respect to $\tau_c$. Therefore, the
$\tau_c$-generator $(\widehat{A},\widehat{D})$ of $(T(t))_{t\geq 0}$
can be defined through the Laplace transform of the semigroup
following the approach of \cite{Kuhnemund03}.

The connection between $A$ and $\widehat A$ is the following.
By \cite[Proposition 3.5, Sections 4,5]{MPW2002},
 it holds that $D_d(A)\subseteq \widehat D \subseteq D_{\max}(A)$ and $\widehat{A}u=Au$ for every $u\in \widehat D$, where $D_{\max}(A)$ is the maximal domain of the operator $A$ in $C_b(\R^N)$, that is
\[
D_{\max}(A):=\bigg\{ \,u\in  \bigcap_{1<p<\infty} W^{2,p}_{\rm{ loc}}(\R^N)\cap C_b(\R^N)\mid\ Au\in C_b(\R^N)\,\bigg\},
\]
while $D_d(A)= D_{{\max}}(A)\cap C_0(\R^N)$ is the Dirichlet domain of $A$.
 In \cite[Theorems 3.7 and 3.12]{MPW2002} sufficient conditions are given in order that $\widehat D=D_{\max}(A)$ or $\widehat D=D_d(A)$.

A probability measure $\mu$ defined on the Borel subsets of $\R^N$ is called an invariant measure for the semigroup $(T(t))_{t\geq 0}$ if for every $f\in C_b(\R^N)$ and $t\geq 0$,
\begin{eqnarray*}
\int_{\R^N}T(t)fd\mu=\int_{\R^N}fd\mu .
\end{eqnarray*}
Whenever an invariant measure exists, $T(t){\bf 1}=\bf 1$ for every
$t\geq 0$ and, hence, $\lambda -A$ is injective on $D_{\max}(A)$ and
the $\tau_c$-generator of $(T(t))_{t\geq 0}$ is $(A,D_{\max}(A))$
(see \cite[Propositions 5.1 and 5.9]{MPW2002}). Since
$(T(t))_{t\geq 0}$ is irreducible and has the strong Feller
property, an invariant measure is unique, if existing (see e.g.,
\cite[Theorem 4.2.1]{dapzab}). Moreover, the invariant measure and
the Lebesgue measure are equivalent (see e.g., \cite[Proposition
8.1.5]{libro}).
In the one dimensional setting, a necessary and sufficient condition
for the existence of an invariant measure is known in the case when
the diffusion coefficient is bounded from below by a positive
constant (see \cite[Proposition 6.2]{MPW2002}. Things are different
in the $N$-dimensional setting. In the particular case when $A$ is an
Ornstein-Uhlenbeck operator, i.e., in the case when
\begin{eqnarray*}
Au(x)=\sum_{i,j=1}^Nq_{ij}D_{ij}u(x)+\sum_{i,j=1}^Nb_{ij}x_jD_iu(x),
\end{eqnarray*}
where $B=(b_{ij})$ and $Q=(q_{ij})$ are $N\times N$ matrices, $Q$
being positively defined, a necessary and sufficient condition for
the existence of an invariant measure of the associated semigroup
(the so called Ornstein-Uhlenbeck semigroup) is available. More
precisely, an invariant measure exists (and in this case is unique)
if and only if the spectrum of $B$ is contained in the left half-plane
$\{\lambda\in\mathbb C: {\rm Re}\,\lambda<0\}$  (see e.g.,
\cite[Section 11.2.3]{DZ}). For more general elliptic operators $A$,
a sufficient condition for the existence of an invariant measure is
the Has'minskii theorem, which requires the existence of a Lyapunov
function $V\in C^2(\R^N)$ blowing up as $|x|\to\infty$, such that
$AV(x)\to -\infty$ as $|x|\to \infty$. For more information about
invariant measures we refer to \cite{dapzab} and \cite{Ro} (see also
\cite{libro} and \cite{Dap01}).

If the semigroup $(T(t))_{t\geq 0}$ has an invariant measure $\mu$,
then it extends to a strongly continuous semigroup of positive
contractions on $L^p(\mu)=L^p(\R^N, d\mu)$ for every $1\leq
p<\infty$, that we still denote by $\Tt$. Moreover,  the measure
$\mu$ is  ergodic (see \cite[Theorem 7.2.2 and Proposition
7.2.1]{Dap01}), i.e., for every $f \in L^2(\mu)$
\begin{equation} \label{eq:80}
 \lim_{r\to \infty} C(r) f = \ov f
\end{equation}
in $L^2(\mu)$, where $\ov f = \int_{\R^N} f \,d\mu$.

Property (\ref{eq:80}) can be extended to $L^p(\mu)$ for every
$1\le p<\infty$ and $f \in L^p(\mu)$. Indeed, if $p>2$, we can
estimate
\begin{align*}
 \|C(r)f -\ov f\|^p_{L^p(\mu)} & = \int_{\R^N} |C(r)f(x) - \ov f|^{2} \cdot |C(r)f(x) - \ov f|^{p-2} \, d\mu(x) \nonumber \\
 &\le \|C(r)f - \ov f\|_{\infty}^{p-2} \cdot \|C(r)f - \ov f\|^2_{L^2(\R^N, \mu)} \nonumber \\
 &\le ( 2 \|f\|_\infty)^{p-2} \cdot \|C(r)f - \ov f\|^2_{L^2(\mu)}.
\end{align*}
On the other hand, if $1\le p <2$, by the H\"older inequality with conjugate exponents $2/p$ and $2/(2-p)$, we have
\begin{align*}
 \|C(r)f - \ov f\|^p_{L^p( \mu)} & = \int_{\R^N} |C(r)f(x) - \ov f|^p \, d\mu(x)  \nonumber \\
 &\le \left ( \int_{\R^N} |C(r)f(x) - \ov f|^{2} \right
)^{\frac{p}{2}}\cdot \mu(\R^N)^{\frac{2-p}{2}}  \nonumber \\ & =
\|C(r)f - \ov f\|_{L^2(\mu)}^{p}.
\end{align*}
As we have already remarked,
\begin{eqnarray*}
\lim_{r\to \infty}C(r)f=\ov f \mbox{ in } L^p(\mu)~\forall f\in
L^p(\mu) \Longleftrightarrow \lim_{\lambda \to 0^+} \lambda
R(\lambda, A)f = \ov f \mbox{ in } L^p(\mu)~\forall f\in L^p(\mu).
\end{eqnarray*}

We are now ready to state and prove the main result of this section.
\begin{theorem}
\label{prop:4.1} The semigroup $\Tt$ is $\tau_c$-mean ergodic if and
only if it admits an invariant measure $\mu$. In this case,
\begin{equation}
 \tau_c-\lim_{r\to \infty}C(r) f =  \ov f = \int_{\R^N} f \, d\mu,
\label{last}
\end{equation}
for every $f \in C_b(\R^N)$.
\end{theorem}

\begin{proof}
To begin with, let us assume that $\Tt$ admits an invariant measure
$\mu$ and let us show that it is $\tau_c$-mean ergodic and satisfies
\eqref{last}. For this purpose, fix $f \in C_b(\R^N)$ and, for each
$\lambda
>0$, set $u_\lambda = \lambda R(\lambda, A)f$. In order to prove the
$\tau_c$-mean ergodicity of the semigroup $\Tt$, by Proposition
\ref{prop:15} it suffices to show that $u_\lambda$ converges with
respect to $\tau_c$ to some function $v \in C_b(\R^N)$  as $\lambda
\to 0^+$.

By the remarks before the statement of the theorem, each function
$u_\lambda$ belongs to $\Dmax (A)$, it solves the equation $\lambda
u_\lambda - A u_\lambda = \lambda f$ and, by estimate \eqref{star},
it satisfies the following estimates
\begin{equation} \label{eq:91}
  \|u_\lambda\|_\infty \le \|f\|_\infty, \qquad \|Au_\lambda\|_\infty \le 2\lambda \|f\|_\infty .
\end{equation}
Moreover, since $u_\lambda$ belongs to $W^{2,p}_{\loc}(\R^N)$ for any $p<\infty$, the following estimate
\begin{equation} \label{eq:92}
  \|u_\lambda\|_{W^{2,p}(B(0,R))} \le C [ \|Au_\lambda\|_{L^p(B(0,2R))} +\|u_\lambda\|_{L^p(B(0,2R))}]
\end{equation}
holds for every $1<p<\infty$, any $R>0$ and some positive constant
$C$ depending only  on $p,R$ and the operator $A$ (see e.g.,
\cite[Theorem 9.11]{gilbtru}).
 Combining (\ref{eq:91}) and (\ref{eq:92}),  we obtain that there exist $\lambda_0>0$ and a constant $C_1>0$ such that
\begin{equation} \label{eq:97}
\|u_\lambda\|_{W^{2,p}(B(0,R))} \le C_1 \|f\|_\infty,
\end{equation}
for all $\lambda \in ]0,\lambda_0]$. It follows that the family
$(u_\lambda)_{\lambda\in ]0,\lambda_0]}$ is uniformly bounded in $W^{2,p}(B(0,R))$ for every
$1<p<\infty$ and every $R>0$. Now, we fix $p>N$ and consider a
sequence $(\lambda_n)_{n \in \N} \subset\, ]0,\lambda_0]$ such that
$\lambda_n \to 0^+$ as $n \to \infty$. By (\ref{eq:97}) and the
Ascoli Arzel\`a theorem, we deduce the existence of a subsequence
$(\lambda_{n_k})_{k\in\N}$ of $(\lambda_n)_{n\in\N}$ and a function
$v \in C_b(\R^N)$ such that $u_{\lambda_{n_k}} \to v$ uniformly on
$\overline{B(0,R)}$ for any $R>0$, as $k \to \infty$. It easily
follows that $(u_{\lambda_{n_k}})_{k\in\N}$ converges to $v$ in
$L^p(B(0,R),\mu)$, too. Since $C(r)f$ tends to $\ov f$ in $L^p(\mu)$
for any $f\in L^p(\mu)$, it follows that $v=\ov f$ $\mu$-a.e. in
$B(0,R)$ and, hence, a.e. in $B(0,R)$ with respect to the Lebesgue
measure as $\mu$ and the Lebesgue measure are equivalent. This is
enough to conclude that $u_\lambda \to  \ov f$ uniformly on compact
sets as $\lambda \to 0^+$. Indeed, if, by contradiction, this were
not true, there would exist a sequence $(\lambda_n)_{n \in \N} \in\,
]0, \lambda_0]$ and two positive constants $R$ and $M$ such that
$\|u_{\lambda_n} - \ov f\|_{C(\overline{B(0,R)})}
> M$ for every $n \in \N$. But this is a contradiction, since the
above arguments applied to the sequence $(u_{\lambda_n})_{n\in\N}$,
show that we can extract a subsequence $(u_{\lambda_{n_k}})_{k\in\N}$
which converges uniformly to $\ov f$ on $B(0,R)$ as $k\to \infty$.

Let us now suppose that $\Tt$ is $\tau_c$-mean ergodic and show that
is admits an invariant measure. Let $P$ be the $\tau_c$-mean ergodic
projection associated with the semigroup $\Tt$. We claim that there
exists $x\in\R^N$ such that $(P\one)(x)\neq 0$. On the contrary,
suppose that this is not the case. Formula \eqref{repres} implies
that $T(t)f\le \|f\|_{\infty}T(t)\one$ for any $f\in C_b(\R^N)$. It
follows that $Pf\le \|f\|_{\infty}P\one$. Hence, if $P\one\equiv 0$,
then $Pf\equiv 0$ as well, for any nonnegative $f\in C_b(\R^N)$.
Splitting a general function $f\in C_b(\R^N)$ into its positive and
negative parts, one easily realizes that $Pf\equiv 0$ for any $f\in
C_b(\R^N)$. Since, by Lemma \ref{lemma-3.5}(i), the range of $P$
equals the kernel of the $\tau_c$-generator of the semigroup $\Tt$,
this would imply that ${\rm ker}\,A=\{0\}$, which clearly is not the
case since the kernel of $A$ contains all the constants.

So, let us fix $x\in\R^N$ such that $P\one(x)\neq 0$. Since $P$ is a
bounded operator in $C_b(\R^N)$ and $\Tt$ is a positivity preserving
semigroup, the map $f\mapsto Pf(x)$ is a positive functional on
$(C_b(\R^N),\|\cdot\|_{\infty})$. By the Riesz representation
theorem (see e.g., \cite[Theorem 2.14]{rud87}), there exists a
finite positive Borel measure $\mu_x$ on $\R^N$ such that
\begin{eqnarray*}
Pf(x)=\int_{\R^N}fd\mu_x,\qquad\;\,f\in C_b(\R^N).
\end{eqnarray*}
Note that $0<\mu_x(\R^N)\le 1$. Indeed, $\mu_x(\R^N)=P\one(x)\in (0,1]$.  Taking advantage of \eqref{eq:65} we can write
\begin{equation}
\int_{\R^N}fd\mu_x=Pf(x)=PT(t)f(x)=\int_{\R^N}T(t)fd\mu_x,
\label{invariant}
\end{equation}
for any $f\in C_b(\R^N)$ and any $t>0$. Up to a normalization, we
can assume that $\mu_x$ is a probability measure. Formula
\eqref{invariant} thus shows that $\mu_x$ is an invariant measure of
$\Tt$.
\end{proof}

\begin{remark}
{\rm
We remark that in \cite[Theorem 4.2.1]{DZ} the authors show that, whenever $\Tt$ admits an invariant measure $\mu$, the
function $T(t)f$ converges to $Pf$ pointwisely in $\R^N$ as $t\to \infty$ for any $f\in C_b(\R^N)$ and hence,  $\Tt$ is mean ergodic with respect to the topology $\tau_p$ of the bounded pointwise convergence.

From all the previous results, we can now show that $\Tt$ is $\tau_c$-mean ergodic
if and only if it is $\tau_p$-mean ergodic.
Clearly, since the topology $\tau_p$ is coarser than the topology $\tau_c$, if $\Tt$ is $\tau_c$-mean ergodic, then it is
$\tau_p$-mean ergodic as well. On the other hand, it is well known that, for any $t>0$, $T(t)$ is continuous from $(X,\tau_p)$ into
$(X,\tau_c)$, i.e., for any bounded sequence $(f_n)_{n\in\N}\in C_b(\R^N)$ which converges pointwisely in $\R^N$ to some function $f\in C_b(\R^N)$, the sequence $(T(t)f_n)_{n\in\N}$ is bounded and converges to $T(t)f$ locally uniformly in $\R^N$ (see
\cite[Proposition 4.6]{MPW2002}). Since the topology $\tau_p$ satisfies Assumptions \ref{ass:1}, if
$\Tt$ is $\tau_p$-mean ergodic, then by Lemma \ref{lemma-3.5} the corresponding $\tau_p$-mean ergodic projection $P_{\tau}$
commutes with $T(t)$ and $T(t)P_{\tau}=P_{\tau}$ for any $t>0$. Now, the arguments in the second part of the proof of Theorem \ref{prop:4.1} show that $\Tt$ admits an invariant measure and, hence, it
is $\tau_c$-mean ergodic. Clearly, the $\tau_p$- and $\tau_c$-mean ergodic projections coincide. In particular, the $\tau_p$-mean (or $\tau_c$-mean)
 ergodicity is then equivalent to the $\tau_p$-convergence of $T(t)$ to a
projection on $C_b(\R^N)$ as  $t\to \infty$.
}
\end{remark}

\begin{remark}
{\rm An example of a not mean-ergodic bi-continuous semigroup  is given by  the left translation (semi)group $(T_l(t))_{t\in\R}$ on $C_b(\R)$. By \cite{K-2,Kuhnemund03} we know  that
$(T_l(t))_{t\in\R}$ is a bi--continuous contraction (semi)group on $(C_b(\R),\|\cdot\|_\infty)$ with respect to $\tau_c$, with  $\tau_c$--generator $(A,D(A))$  given  by
$Au=u'$, for $f\in D(A)=C^1_b(\R)$ (
$C^1_b(\R)$ being the space of all the differentiable  functions on $\R$, with continuous and bounded derivatives). In particular,  $(T_l(t))_{t\in\R}$ is not be $\tau_c$--mean ergodic. Indeed,
 a
counterexample is given in \cite[Chapter V, Examples 4.12]{Nagel}.
Let $f\in C_b(\R)$ be defined piecewisely as follows:
\begin{eqnarray*}
f(x)= \left\{
\begin{array}{ll}
2(-1)^k(x-10^k+1/2), & x\in [10^k-1/2,10^k[\,,\\[1mm]
(-1)^k, & x\in [10^k,10^{k+1}-1[\,,\\[1mm]
2(-1)^{k+1}(x-10^{k+1}+1/2), & x\in [10^{k+1}-1,10^{k+1}-1/2],\\[1mm]
0, & \mbox{elsewhere in }\R.
\end{array}
\right.
\end{eqnarray*}
A straightforward computation reveals that
\begin{align*}
\frac{1}{10^{n+1}-1/2}\int_0^{10^{n+1}-1/2}T(t)f(0)dt&=\frac{1}{10^{n+1}-1/2}\int_0^{10^{n+1}-1/2}f(t)dt\\
&=(-1)^{n+1}\frac{9}{11}\frac{10^{n+1}}{10^{n+1}-1/2}+o(1),
\end{align*}
where $o(1)$ denotes terms which tend to $0$ as $n\to\infty$. It
follows that $C(r)f(0)$ does not converge as $r\to \infty$,
implying that the semigroup is not $\tau_c$-mean ergodic.

Another proof of the fact that $\big( T_l(t) \big)_{t\geq 0}$ is not
$\tau_c$-mean ergodic can be obtained adapting the arguments in the
second part of the proof of Theorem \ref{prop:4.1}. Indeed, by
contradiction suppose that $\big( T_l(t) \big)_{t\geq 0}$ is
$\tau_c$-mean erdodic. Then, the projection $P$ can be identified
with a positive functional on $C_b(\R)$ since the kernel of the
operator $Au=u'$ consists of constant functions only. Therefore,
there exists a finite positive Borel measure $\mu$ such that
\begin{equation}
\int_{\R}f(t+\cdot)d\mu=\int_{\R}fd\mu,\qquad\;\,t>0,
\label{ultima-ultima}
\end{equation}
for any $f\in C_b(\R)$. By density, the previous formula can be
extended to any bounded Borel measurable function $f$. Hence,
writing \eqref{ultima-ultima} with $t=n\in\mathbb N\cup\{0\}$ and
$f=\chi_{[0,1[}\,$, we obtain that $\mu([n,n+1[)=\mu_{[0,1[}\,$. It
follows that
\begin{eqnarray*}
\infty>\mu([0,\infty[)=\sum_{n=0}^{\infty}\mu([n,n+1[)=\sum_{n=0}^{\infty}\mu([0,1[),
\end{eqnarray*}
implying that $\mu([0,\infty[)=0$. Applying the same
argument to the function $f_{-k}=\chi_{[-k,-k+1[}$, we get
$\mu([-k,-k+1[)=0$ for any $k\in\N$, so that $\mu(]-\infty,0[)=0$. As a byproduct, it follows that
$0=\mu(\R)=P\one$. But this is a contradiction since it implies that
$Pf=0$ for any $f\in C_b(\R)$, i.e., ${\rm ker}A=\{0\}$, which is
not the case.}
\end{remark}

\subsection{A second application: the nonautonomous case}

Let now consider the case when the coefficients of the operator $A$ depend on time, i.e., the case when
\begin{eqnarray*}
A(s)u(x) = \sumij q_{ij}(s,x)D_{ij} u(x) + \sum_{i=1}^N b_i(s,x) D_i
u(x), \qquad s\in\R,\;\,x \in \R^N.
\end{eqnarray*}
In \cite{KLL,LLZ,LZ} such operators have been extensively studied both in periodic and nonperiodic settings.
Here, we confine ourselves to the case of $T$-periodic coefficients analyzed in \cite{KLL,LZ}. We denote by ${\mathbb T}$
the torus $\mathbb T=[0,T]$ {\it mod.} $T$, and
assume the following standing hypotheses on the coefficients:
\begin{enumerate}[\rm 1)]
\item
the coefficients $q_{ij}$ and $b_i$ $(i,j=1,\ldots,N)$ are $T$-time
periodic and belong to $C^{\alpha/2,\alpha}_{\rm loc}(\mathbb T
\times
 \R^N)$ for any $i,j=1,\ldots,N$
and some $\alpha\in (0,1)$;
\item
for every $(s,x)\in  \R^{1+N}$, the matrix $Q(s,x)=(q_{ij}(s,x))$ is
symmetric and there exists a function $\eta: {\mathbb
T}\times\R^N\to \R$ such that $0<\eta_0:=\inf_{\mathbb
T\times\R^N}\eta$ and
\begin{eqnarray*}
\langle Q(s,x)\xi , \xi \rangle \geq \eta(s,x) |\xi |^2, \qquad\;\,
\xi \in \R^N,\;\,(s,x) \in \mathbb T \times \R^N;
\end{eqnarray*}
\item
there exist a positive function $V  \in C^2(\R^N)$ and constants
$a$,  $c>0$    such that
\begin{eqnarray*}
\qquad\;\, \lim_{|x|\to \infty} V (x) = \infty \quad \mbox{
and } \quad (A(s) V)( x) \leq a - c V(x), \quad (s,x)\in  \mathbb T \times \R^N.
\end{eqnarray*}
\end{enumerate}

Under the previous set of assumptions, in \cite{KLL} it has been shown that an evolution operator
$(P(t,s))$ can be associated with the operator $A(t)$ in $C_b(\R^N)$. For any $f\in C_b(\R^N)$, $P(r,s)f$ is the value at $r$ of
the unique bounded classical solution to the Cauchy problem
\begin{eqnarray*}
\left\{
\begin{array}{ll}
D_tu(t,x)=A(t)u(t,x), & (t,x)\in\, ]s,\infty[\times\R^N,\\[1mm]
u(s,x)=f(x), &x\in\R^N.
\end{array}
\right.
\end{eqnarray*}
Here, by classical solution we mean a function $u\in
C^{1,2}(]s,\infty[\times\R^N)\cap C([s,\infty[\,\times\R^N)$ which
solves the above Cauchy problem. A variant of the classical maximum
principle, which can be proved under assumption 3), shows that
$P(t,s)$ is a contraction in $C_b(\R^N)$ for any $s\le t$. Moreover,
since the coefficients of the operator $A(t)$ are $T$-periodic,
$P(t+T,s+T)f=P(t,s)$ for any $s,t\in\R$ such that $s<t$.

Further, the previous set of assumptions imply that a periodic evolution system of measures (also
called {\it entrance laws at $-\infty$} in \cite{dynkin}) can be associated with the evolution operator $P(t,s)$, i.e., there exists a one-parameter
family of probability measures $\{\mu_s: s\in\R\}$ such that
\begin{eqnarray*}
\int_{\R^N}P(t,s)fd\mu_t=\int_{\R^N}fd\mu_s,\qquad\;\,s<t,
\end{eqnarray*}
for any $f\in C_b(\R^N)$. Evolution systems of measures are the natural counterpart of invariant measures for nonautonomous operators.

By \cite[Proposition 2.10]{LLZ}, the previous one is the only periodic evolution system of measures of $(P(t,s))$. Note that
the nonperiodic evolution systems of measures are, in general, infinite many. See for instance \cite{GL} where
they have been all characterized in the particular case when $A(t)$ is the nonautonomous Ornstein-Uhlenbeck operator.

Let us define a measure in $\mathbb T\times\R^N$ by extending to the $\sigma$-algebra of all the Borel sets
of $\mathbb T\times\R^N$ the set-valued function defined by
$$
\mu(A\times B)=\frac{1}{T}\int_A\mu_s(B)ds,
$$
for any pair of Borel sets $A\subset\mathbb T$ and $B\subset\R^N$.
The measure $\mu$ turns out to be the (unique) invariant measure of the evolution semigroup (also called Howland semigroup) $({\mathcal T}(t))_{t\geq 0}$ defined by
\begin{eqnarray*}
(\mathcal T(t)f)(s,x)=(P(s,s-t)f(s-t,\cdot))(x),\qquad\;\,(s,x)\in\R^{1+N},
\end{eqnarray*}
for any $f\in C_b(\mathbb T\times\R^N)$. It follows immediately that
$({\mathcal T}(t))_{t\geq 0}$ can be extended to a strongly
continuous semigroup to $L^p(\mathbb T\times\R^N,\mu)$ for any $p\in
[1,\infty[$\,. Moreover, $(\mathcal T(t))_{t\ge 0}$ is a
bi-continuous semigroup on $C_b(\mathbb T\times\R^N)$ with respect
to the topology of uniform convergence on compact sets of $\mathbb
T\times\R^N$, each operator ${\mathcal T}(t)$ being a
$\|\cdot\|_\infty$-contraction.

If we denote by $G_p$ the infinitesimal generator of the semigroup $(\mathcal T(t))_{t\ge 0}$ in $L^p(\R^{1+N},\mu)$, it follows
that $G_p$ extends the operator
${\mathcal G}\varphi:=A(t)\varphi-D_t\varphi$ defined on smooth and $T$-periodic (with respect to $t$) functions $\varphi$.
In particular, \cite[Proposition 6.7]{LZ} shows that the set
\begin{eqnarray*}
D(G_{\infty})  = \bigg\{f\in\bigcap_{q<\infty} W^{1,2}_{q,\rm loc}( \mathbb T \times \R^N, ds\times dx)
\cap C_b(\mathbb T \times\R^N),\;{\mathcal G} f\in C_b(\mathbb T\times \R^N)\bigg\},
\end{eqnarray*}
is a core of $G_p$. $D(G_{\infty})$ is nothing but the domain of the weak generator of the restriction of $(\mathcal T(t))_{t\ge 0}$ to
$C_b(\mathbb T\times\R^N)$.

The proof of \cite[Proposition 2.10]{LLZ} (which is based on \cite[Proposition 3.2.5]{DD}) shows
that ${\mathcal T}(t)\chi_{\Gamma}=\chi_{\Gamma}$ for any $t>0$ if and only if the Borel set $\Gamma\subset\mathbb T\times\R^N$
satisfies either $\mu(\Gamma)=0$ or $\mu(\Gamma)=1$.
Using this property and arguing as in the autonomous case (see e.g., \cite[Proposition 8.1.11]{libro}) one can easily show
that $\mu$ is ergodic, i.e.,
\[
\lim_{r\to\infty}\frac{1}{t}\int_0^t{\mathcal T}(s)ds=\int_{\mathbb
T\times\R^N}fd\mu,
\]
in $L^2({\mathbb T}\times\R^N,\mu)$ for any $f\in L^2({\mathbb
T}\times\R^N,\mu)$. The same arguments as in the previous subsection
allow us to prove that the integral average of ${\mathcal
T}(\cdot)f$ to $\overline f$ converges to $\overline f$ in
$L^p({\mathbb T}\times\R^N,\mu)$ for any $f\in L^p({\mathbb
T}\times\R^N,\mu)$ and any $p\in [1,\infty[$\,.

Since no confusion may arise, we denote indifferently by $\tau_c$ both the topology of uniform convergence on compact sets of  $\mathbb T\times\R^N$ and the topology of uniform convergence on compact sets of $\R^N$.
We can now prove the following result.

\begin{proposition}
Under the previous set of assumptions, the semigroup $(\mathcal
T(t))_{t\ge 0}$ is $\tau_c$-mean ergodic. In particular,
\begin{equation}
\tau_c\mbox{\rm --}\lim_{t\to
\infty}\frac{1}{t}\int_0^t(P(s,-r)f-m_{-r}f)dr=0,
\label{stima-media}
\end{equation}
for any $r,s>0$. Moreover, a Lyapunov type theorem holds. More
specifically, the constants are the only solutions to the parabolic
equations $D_tu-A(t)u=0$ which belong to $\bigcap_{q<\infty}
W^{1,2}_{q,\rm loc}( \mathbb T \times \R^N, ds\times dx) \cap
C_b(\mathbb T \times\R^N)$.
\end{proposition}

 \begin{proof}
 Following essentially the same lines as in the proof of Proposition \ref{prop:4.1}
 we can show that the function $u_{\lambda}=\lambda R(\lambda,G_{\infty})f$
converges to $\overline f=\int_{\mathbb T\times\R^N}fd\mu$ locally uniformly in $\R^{1+N}$ as $\lambda\to 0^+$.
Indeed, the arguments in the quoted proof shows that
$\|{\mathcal G}u_{\lambda}\|_{\infty}+\|u_{\lambda}\|_{\infty}\le 3\|f\|_{\infty}$ for any $\lambda\in (0,1]$.
By the parabolic $L^p$-interior estimates, for any $R>0$, there exists a positive constant $C=C(p,R)$ such that
\begin{align*}
\|u_{\lambda}\|_{W^{1,2}_p((0,T)\times B(0,R))}&\le C\left (\|u\|_{L^p((0,T)\times B(0,2R))}+\|{\mathcal G}u\|_{L^p((0,T)\times B(0,2R))}\right )\\
&\le C'\left (\|u\|_{\infty}+\|{\mathcal G}u\|_{\infty}\right )
\le C''\|f\|_{\infty},
\end{align*}
for any $p\in\, ]1,\infty[$\,. Taking $p$ sufficiently large, we
can conclude that the family of functions $u_{\lambda}$ ($\lambda\in
(0,1]$) is bounded in $C^{\alpha}([0,T]\times B(0,R))$ for any
$R>0$. From now on the proof can be carried over as in the proof of
Proposition \ref{prop:4.1} taking into account that the measure
$\mu$ is equivalent to the Lebesgue measure (see \cite[Theorem
3.8]{BKR}). This proves the first part of the proposition.

To complete the proof, we observe that, since $C_b(\R^N)$
continuously embeds into $C_b(\mathbb T\times\R^N)$, and $(\mathcal
T(t)f)(s,x)=(P(s,s-t)f)(x)$ for any $(s,x)\in\R^{1+N}$ and any
$t>0$, the result so far proved shows that
\begin{eqnarray*}
\tau_c\mbox{--}\lim_{t\to\infty}\frac{1}{t}\int_0^tP(s,s-r)fdr=\frac{1}{T}\int_0^Tds\int_{\R^N}fd\mu_s.
\end{eqnarray*}

Let us now observe that
\begin{equation}
\frac{1}{T}\int_0^Tds\int_{\R^N}fd\mu_s=\lim_{t\to\infty}\frac{1}{t}\int_0^tds\int_{\R^N}fd\mu_s.
\label{limit}
\end{equation}
Of course, it suffices to prove the previous property for nonnegative functions. The general case then will follow splitting $f$
into its positive and negative parts. Fix $t>T$ and split $t=nT+r$ for some $n\in\N$ and some $r\in [0,T)$.
Then,
\begin{eqnarray*}
\frac{1}{t}\int_0^tds\int_{\R^N}fd\mu_s\ge \frac{1}{(n+1)T}\int_0^{nT}ds\int_{\R^N}fd\mu_s
=\frac{n}{(n+1)T}\int_0^Tds\int_{\R^N}fd\mu_s,
\end{eqnarray*}
where we have used the $T$-periodicity of the $T$-valued function $s\mapsto\mu_s$. Hence,
\begin{eqnarray*}
\liminf_{t\to \infty}\frac{1}{t}\int_0^tds\int_{\R^N}fd\mu_s
\ge\frac{1}{T}\int_0^Tds\int_{\R^N}fd\mu_s.
\end{eqnarray*}
On the other hand,
\begin{eqnarray*}
\frac{1}{t}\int_0^tds\int_{\R^N}fd\mu_s\le \frac{1}{nT}\int_0^{(n+1)T}ds\int_{\R^N}fd\mu_s
=\frac{n+1}{nT}\int_0^Tds\int_{\R^N}fd\mu_s,
\end{eqnarray*}
which shows that
\begin{eqnarray*}
\limsup_{t\to \infty}\frac{1}{t}\int_0^tds\int_{\R^N}fd\mu_s
\ge\frac{1}{T}\int_0^Tds\int_{\R^N}fd\mu_s.
\end{eqnarray*}
Formula \eqref{limit} now follows. Similarly, a straightforward change of variables shows that
\begin{align*}
\frac{1}{t}\int_0^tP(s,s-r)fdr&=\frac{1}{t}\int_{-s}^{t-s}P(s,-r)fdr\\
&=\frac{1}{t}\int_0^tP(s,-r)fdr+\frac{1}{t}\int_{-s}^0P(s,-r)fdr
-\frac{1}{t}\int_{t-s}^tP(s,-r)fdr,
\end{align*}
for any $s\ge 0$, and the last two integral terms in the previous chain of equalities tends to $0$ as $t\to \infty$ uniformly in
$[0,T]\times\R^N$. This shows that
\begin{eqnarray*}
\tau_c\mbox{--}\lim_{t\to \infty}\frac{1}{t}\int_0^tP(s,s-r)fdr=\frac{1}{t}\int_{-s}^{t-s}P(s,-r)fdr
=\tau_c\mbox{--}\lim_{t\to \infty}\frac{1}{t}\int_0^tP(s,-r)fdr,
\end{eqnarray*}
accomplishing the proof of \eqref{stima-media}. The other statement of the proposition follows from Lemma 2.5(i).
\end{proof}

\begin{remark}\rm
Formula \eqref{stima-media} states a convergence to $0$ in integral
mean of $P(t,s)f-m_sf$. Under stronger assumptions than those we are
assuming here, it is possible to prove that
\begin{eqnarray*}
\tau_c\mbox{--}\lim_{s\to\infty}(P(t,s)f-m_sf)=0,
\end{eqnarray*}
for any fixed $t$. We refer the reader to \cite{LLZ} for more
details.
\end{remark}

\end{document}